\documentclass[12pt]{article}
\title{A generator system of invariant differential forms.}
\author{{\rm Tensai Bakabon}\\
{\tt\small Bakada University, Department of Mathematics}\\
{\tt\small 3-10-5 Tsukimidai Susukigahara, Nerima-ku, Tokyo, Japan}\\
{\tt\small e-mail:korede\_ii\_noda@yahoo.co.jp}}
\usepackage{amssymb}
\date{April 1,2002}

\newtheorem{theorem}{Theorem}

\begin{document}

\maketitle

\begin{abstract}
We obtain a generator system of the algebra of $\mathrm{GL}(V)$-invariant 
differential forms on $\mathrm{End} _{\bf k} (V)$. 
The proof uses the Weyl-Schur reciprocity.
\end{abstract}

Let ${\bf k}$ be a field of characteristic zero and 
$\mathrm{M}_n({\bf k})$ be the set of all $n\times n$ matrices with 
entries of elements of ${\bf k}$. We consider $\mathrm{GL}_n({\bf k})$ adjoint 
action on $\mathrm{M}_n({\bf k})$, i.e. for given $g\in \mathrm{GL}_n({\bf k})$ 
and $X\in \mathrm{M}_n({\bf k})$, $\mathrm{ad}(g)X:=g\cdot X\cdot g^{-1}$.
It is classically known that the algebra 
$\mathbb{C}[\mathrm{M}_n({\bf k})]^{\mathrm{GL}_n({\bf k})}$
of all 
$\mathrm{GL}_n({\bf k})$-invariant polynomials on $\mathrm{M}_n({\bf k})$ is 
generated by $\mathrm{Tr}(X^\ell )\ (\ell=1,2,\cdots )$, where 
$X=(x_{i j})_{i,j=1}^{n}$. More precisely, 
if we define polynomials $\sigma_1(X),\cdots ,\sigma_n(X)$ of $X$ by
\begin{displaymath}
\det(t\cdot I_n +X)=t^n +\sigma_1(X)t^{n-1} +\cdots +\sigma_n(X),
\end{displaymath}
the invariant ring is the polynomial ring generated by
$\sigma_1(X),\cdots ,\sigma_n(X)$ over ${\bf k}$. Furthermore 
polynomials $\sigma_1(X),\cdots ,\sigma_n(X)$ and 
$\mathrm{Tr}(X^\ell )\ (\ell=1,2,\cdots )$ are connected by Newton's formula.

We next consider $\mathrm{GL}_n({\bf k})$-invariant 
differential forms on $\mathrm{M}_n({\bf k})$ with polynomial coefficients. 
Similar to the above case, we obtain many invariant 
differential forms of degree $p$ by the formula 
\begin{eqnarray*}
\omega^{\ell_1\ell_2\cdots \ell_p} &:=& \mathrm{Tr}\left(X^{\ell_1}dX\wedge X^{\ell_2} dX\wedge\cdots \wedge X^{\ell_p} dX \right) \\
&=& \sum_{i_1,\cdots ,i_p, j_1,\cdots ,j_p=1}^{n}(X^{\ell_1})_{i_1 j_1}dx_{j_1 i_2}\wedge (X^{\ell_2})_{i_2 j_2}dx_{j_2 i_3}\cdots (X^{\ell_p})_{i_p j_p}dx_{j_p i_1},
\end{eqnarray*}
where $(X^{\ell})_{i j}$ is the $(i,j)$
component of $X^\ell$ and $dX=(dx_{ij})_{i,j=1}^{n}$.

It is natural to ask if the algebra of all $\mathrm{GL}_n({\bf k})$-invariant 
differential forms is generated by these forms.
The purpose of this paper is to answer this question affirmatively.

\begin{theorem}
The algebra of all $\mathrm{GL}_n({\bf k})$-invariant differential forms 
on $\mathrm{M}_n({\bf k})$ is generated by 
$\omega^{\ell_1\ell_2\cdots \ell_p}=\mathrm{Tr}(X^{\ell_1}dX\wedge X^{\ell_2} \cdots \wedge X^{\ell_p} dX),\ (p\geq 0,\ 
\ell_1,\cdots ,\ell_p \geq 0)$ over ${\bf k}$. In other words, 
any invariant form is expressed as a linear combination of 
forms of the above type and their product.
\end{theorem}
{\bf Proof of Theorem.}\\
Put $V={\bf k}^n$, $E:=\mathrm{End}_{\bf k}(V)\cong \mathrm{M}_n({\bf k})
\cong V\otimes V^\ast$ 
and $G:=\mathrm{GL}(V)$.
Then we can identify $E$ and its dual $E^\ast$ because 
$E\ni X,Y\mapsto \mathrm{Tr}(XY)$ is a $\mathrm{GL}(V)$-invariant 
non-degenerate bilinear pairing. So the set of all differential forms 
on $E$ is 
\begin{displaymath}
\Omega (E)\cong \mathcal{S}(E)\otimes \Lambda (E)
\cong \bigoplus_{p,q\geq 0}\mathcal{S}^p(E)\otimes \Lambda ^q(E)
\end{displaymath}
where $\mathcal{S}(E)$ and $\Lambda (E)$, respectively, are symmetric and 
alternative tensor algebra of the vector space $E$. 
Hence it is sufficient to show that each component 
$\mathcal{S}^p(E)\otimes \Lambda ^q(E)$ is generated by 
$\omega^{\ell_1\ell_2\cdots \ell_r}\ $
$(\ell_1 +\cdots +\ell _r\leq p,\ r\leq q)$. 
Consider the $G$-equivariant projection map 
$\mathcal{P}:E^{\otimes (p+q)}\rightarrow \mathcal{S}^p(E)\otimes \Lambda ^q(E)$
\begin{displaymath}
\begin{array}{ccl}
E^{\otimes (p+q)}&\ni &v_1 \otimes \cdots \otimes v_p \otimes w_1\otimes \cdots \otimes w_q\\
\downarrow &&\downarrow\\
\mathcal{S}^p(E)\otimes \Lambda ^q(E)&\ni &\frac{1}{p!q!}\sum\limits_{\sigma \in \mathfrak{S}_p,\ \tau \in \mathfrak{S}_q}\mathrm{sgn}(\tau) v_{\sigma ^{-1}(1)}\otimes \cdots \otimes v_{\sigma ^{-1}(p)} \otimes w_{\tau ^{-1}(1)}\otimes \cdots \otimes w_{\tau ^{-1}(q)}.
\end{array}
\end{displaymath}
Since $G$ is a reductive group, $\mathcal{P}$ induces a 
surjection of $G$-invariant tensors to $G$-invariant differential forms.
\begin{displaymath}
\mathcal{P}|_{(E^{\otimes (p+q)})^G}:(E^{\otimes (p+q)})^G 
\twoheadrightarrow (\mathcal{S}^p(E)\otimes \Lambda ^q(E))^G.
\end{displaymath}
Thus we reduce the problem to description of $G$-invariant tensors 
and its images by $\mathcal{P}$. Recall 
$E^{\otimes r} \cong \mathrm{End}_{\bf k}(V^{\otimes r})$ and 
$(E^{\otimes r})^G \cong \mathrm{End}_G(V^{\otimes r})$, 
we have a special element $\Phi (\rho )$ of $G$-invariant 
tensors $(E^{\otimes r})^G$ 
corresponding to a permutation $\rho \in \mathfrak{S}_r$,
\begin{displaymath}
\mathrm{End}_G(V^{\otimes r})\ni \Phi (\rho ):
e_1\otimes \cdots e_r \mapsto e_{\rho ^{-1}(1)}\otimes \cdots e_{\rho ^{-1}(r)}.
\end{displaymath}

Recall that the Weyl-Schur reciprocity assert that 
\begin{displaymath}
\{ \Phi (\rho) \in (E^{\otimes r})^G \ |\ \rho \in \mathfrak{S}_r\}
\mbox{ form a basis of }(E^{\otimes r})^G.
\end{displaymath}

Put $r=p+q$, it is easily checked that 
$\mathcal{P}(\Phi(\rho))\in (\mathcal{S}^p(E)\otimes \Lambda ^q(E))^G$ is a 
product of $\omega^{\ell_1\ell_2\cdots \ell_s}$s for appropriate $s$ and 
$\ell_1,\ell_2,\cdots ,\ell_s$. In particular, the number of factors 
of $\mathcal{P}(\Phi(\rho))$ is equal to the number of cyclic permutations 
in the cyclic decompositon of $\rho \in \mathfrak{S}_r$.
This completes the proof. 
\begin{flushright}
$\square$
\end{flushright}

Perhaps it helps some calculations in relative de Rham cohomologies of 
invariant maps \cite{ve}.

\bigskip

{\large\bf Acknowledgement}\\
The author thanks Ms.Ys.(real) for decisive help in the 
proof of our main result.



\end{document}